\documentclass[twoside, 11pt]{article}
\usepackage{mathrsfs}

\usepackage{mathrsfs,amsfonts,amsmath}
\RequirePackage{ifthen,calc}
\usepackage{theorem}
\usepackage[dvips]{color}

 \catcode`@=11
 \def\@evenhead{\hbox to\textwidth{\footnotesize\rm\thepage \hfill
  {\it }}} % authors name

 \def\@oddhead{\hbox to \textwidth{\footnotesize{\it
  Brownian motion between two random trajectories  } \hfill\thepage}}% abbreviate title

 \renewcommand{\section}{\makeatletter
 \renewcommand{\@seccntformat}[1]{{\csname the##1\endcsname.}\hspace{0.45em}}
 \makeatother \@startsection
{section}%                                            the name
{1}%                                                  the level
{0pt}%                                                the indent
{\baselineskip}%                                      the beforeskip
{0.5\baselineskip}%                                   the afterskip
{\normalsize\bfseries\mathversion{bold}}}
\catcode`@=12
\newcommand\ack{\section*{Acknowledgement}}

\newtheorem{thm}{\noindent Theorem}[section]
\newtheorem{lem}{\noindent Lemma}[section]
\newtheorem{cor}{\noindent Corollary}[section]
\newtheorem{prop}{\noindent Proposition}[section]

\theoremheaderfont{\normalfont\bfseries} \theorembodyfont{\slshape}
{\theorembodyfont{\rmfamily}}
{\theorembodyfont{\rmfamily}}
{\theorembodyfont{\rmfamily}}
{\theorembodyfont{\rmfamily}}
{\theorembodyfont{\rmfamily}}
{\theorembodyfont{\rmfamily}\newtheorem{rem}{\noindent Remark}[section]}
{\theorembodyfont{\rmfamily}}
{\theorembodyfont{\rmfamily}}
{\theorembodyfont{\rmfamily}}
 \def\beqlb{\begin{eqnarray}}\def\eeqlb{\end{eqnarray}}
 \def\beqnn{\begin{eqnarray*}}\def\eeqnn{\end{eqnarray*}}
 \numberwithin{equation}{section}
\def\qed{\hfill$\square$\smallskip}

\def\bfE{{\mathbb{E}}}
\def\bfP{{\mathbb{P}}}
\def\bfR{{\mathbb{R}}}

\def\bfN{{\mathbb{N}}}

\begin{document}

\title{\bf Brownian motion between two random trajectories
}
\author{ You Lv\thanks{Email: youlv@mail.bnu.edu.cn }
\\ \small School of science of mathematics, Beijing Normal University,
\\ \small Beijing 100875, P. R. China.
}
\date{}
\maketitle

%\renewcommand{\thefootnote}{\fnsymbol{footnote}}\footnotetext[1]{}
%\renewcommand{\baselinestretch}{1.0}
%\noindent\hrulefill

\noindent\textbf{Abstract}: Consider the first exit time of one-dimensional Brownian motion $\{B_s\}_{s\geq 0}$ from a random passageway. We discuss a Brownian motion with two time-dependent random boundaries in quenched sense. Let $\{W_s\}_{s\geq 0}$ be an other one-dimensional Brownian motion independent of $\{B_s\}_{s\geq 0}$ and let $\bfP(\cdot|W)$ represent the conditional probability depending on the realization of $\{W_s\}_{s\geq 0}$.
We show that
$$-t^{-1}\ln\bfP^x(\forall_{s\in[0,t]}a+\beta W_s\leq B_s\leq b+\beta W_s|W)$$ converges to a finite positive constant $\gamma(\beta)(b-a)^{-2}$ almost surely and in $L^p~ (p\geq 1)$ if
$a<B_0=x<b$ and $W_0=0.$ When $\beta=1, a+b=2x,$ it is equivalent to the random small ball probability problem in the sense of equiditribution, which has been investigated in \cite{DL2005}. We also find some properties of the function $\gamma(\beta)$.  An important moment estimation has also been obtained, which can
be applied to discuss the small deviation of random walk with random environment in time (see \cite{Lv2018}).
%We investigate the probability of a Brownian motion staying between two trajectories which are related to another Brownian motion.

\smallskip


\smallskip

\noindent\textbf{2000 Mathematics Subjects
Classification:} 60J65

\section{Introduction}
The first exit time of Brownian motion is a classic and interesting topic which has been researched by many scholars.
Let us first recall a very basic result in this field. For a standard Brownian motion $\{B_t\}_{t\geq 0}$ starting from $x$, it is known that
\beqlb\label{sec-1}
\lim\limits_{t\rightarrow+\infty}\frac{-\ln\bfP^x(\forall_{s\leq t}a\leq B_s\leq b)}{t}=\frac{\pi^2}{2(b-a)^2},
\eeqlb
where $a<x<b$. (1.1) shows that the first exit time from a bounded interval has negative exponential tail distribution and the coefficient depends on the width of the bounded interval.

 A lot of further work has been done on the first exit time of Brownian motion. For the Brownian motion in high dimensional space, \cite{LS2002},~\cite{L2003} researched the first exit time $T$ from a fixed convex domain, showing that $-\ln\bfP(T>t)=O(t^{\alpha}),$ where $\alpha$ is a positive constant depending on the degree of dimension and the shape of the convex field. Another extension is to consider the time-dependent boundary. \cite{N1981} studied the asymptotic behavior of
$\bfP(\forall_{0\leq s\leq t}|B_t|\leq f(t)),$ where the boundaries $-f(t)$ and $f(t)$ depend on $t.$ The work \cite{XZ2015} considered the Brownian motion with two linear boundaries and calculated the distribution of the Brownian motion hitting the upper boundary before hitting the lower boundary.
The model of Brownian motion with two time-dependent boundaries can be applied to many different fields such as finance (see \cite{NFK2003}),  biophysical models (see \cite{R1977}) and statistical sequential analysis (see \cite{S1985}).

There are several profound conclusions when the boundaries not only depend on $t$ but also a random variable. \cite{LS2002}, \cite{L2003} and \cite{LSl2012} all discussed the probability
$$\bfP(\forall_{0\leq s\leq t} \parallel B_s\parallel^p\leq 1+ \mu s^r +W(s)),$$
where $\mu>0, r\in[0,1),p>1. ~\{B_s\}_{s\geq 0}$ is a $d$-dimensional ($d\geq 2$) standard Brownian motion and $\{W_s\}_{s\geq 0}$ is a one-dimensional Brownian motion which is independent of $\{B_s\}_{s\geq 0}$. $``\parallel\cdot\parallel"$ is the Euclidean norm. That can be viewed as the first exit time $T$ of Brownian motion $\{B_s\}_{s\geq 0}$ from a random domain. Under some suitable conditions, they all showed that $t^{-\frac{p-1}{p+1}}\ln\bfP(T>t)$ converges to a negative constant which depends on the random domain and on the dimension as $t\rightarrow +\infty.$

What we want to discuss is the decay rate of \beqlb\label{sec-1.2}\bfP^x(\forall_{s\in[0,t]}a+\beta W_s\leq B_s\leq b+\beta W_s|W)\eeqlb as $t\rightarrow +\infty,$ where $\beta\geq 0, a<x<b$ and $\{W_t\}_{t\geq 0}$ is an other standard Brownian motion which is independent of $\{B_t\}_{t\geq 0}.$ Of course, it can also be viewed as the first exit time from a random and time-dependent passageway. We obtain a kind of quenched
result. We prove that the decay rate is $e^{-ct}$ almost surely and $c=(b-a)^{-2}\gamma(\beta).$  Moreover, the function $\gamma(\beta)$ is strictly increasing on $[0,+\infty),$  that is to say, although the width of the random passageway is always constant $``b-a"$, more violent fluctuation of the center will make the first exit time $T$ much shorter. We should notice that if $\beta=1$ and $x=\frac{a+b}{2},$ by scaling property of the Brownian motion, (1.2) has the same distribution as $$\bfP^0(\forall_{s\in[0,1]} |B_s-W_s|\leq \frac{b-a}{2\sqrt{t}}|W).$$ That is the random small ball probability which has been investigated in \cite{DL2005}. So we can see the convergence (2.2) holds in probability from \cite[Theorem 6.1]{DL2005} if $\beta=1$ and $x=\frac{a+b}{2}$.
Motivated by the precise asymptotics of a random quantization problem (see \cite{DFMS2003}), \cite{D2003} and \cite{DL2005} first investigated the random small ball probabilities and gave many important asymptotic estimations for the Gaussian measure $\mu$ on a set centered at a random trajectory when the distribution of the random trajectory is also $\mu.$ In the proof of \cite[Theorem 6.1]{DL2005}, the observing of subadditivity also gives us essential inspiration. Compared with \cite[Theorem 6.1]{DL2005}, our difference is that we also consider the situation of $\beta\neq 1.$ Moreover, we conclude that the convergence (2.2) is also almost surely and uniform for the location of the starting point and the width of the interval at the last moment. By the way, we obtain a moment estimation
(see Theorem 3.1). All of these adjustments will play key role on the research of the small deviation for random walk with random environment in time (see \cite{Lv2018}), which is a main application and motivation of this paper. One can utilize the Brownian motion between two random trajectories to approximate the random walk with random environment in time. Furthermore, the result of \cite{Lv2018} will be a basic tool when we study the barrier problem of the branching random walk with random environment in time. The latter is a work in progress.

Another important point is that our main result can also be viewed as an extension of \cite[Theorem1.1]{MM2015}. Mallein and Mi{\l}o\'{s} consider the probability of a Brownian motion staying above a trajectory of another
 independent Brownian motion. To be more precisely, they proved $$-\ln\bfP(\forall_{0\leq s \leq t} B_s\geq W_s-1|W)/\ln t\rightarrow\gamma, ~\gamma>\frac{1}{2}$$ almost surely and in $L^p ~(p\geq 1)$. The idea of our proof is partly inspired by \cite{MM2015}. However, we face new difficulties
 when we do the moment estimation (see Theorem 3.1) since the probability of Brownian motion with two boundaries is usually smaller than the single
 boundary case.

 The rest of this paper is organized as follows. We state the main theorem and corollaries in section 2. An important estimation of tail distribution is obtained in section 3. Based on this estimation, we give the proof of the main theorem and corollaries in section 4 and section 5 respectively.

\section{Main result}
Throughout this paper, we assume that real numbers $a,b,a',b',a_0,b_0 $ meet the following {\bf basic relationship}
\beqlb\label{sect-0} a<a_0\leq b_0<b,~~ a\leq a'< b'\leq b. \eeqlb

%$$ m,n,k\in \bfN^+\cup\{0\},~~ m\leq n,~~ s,t \in\bfR^+\cup\{0\}. ~~~~~~~~~~~~~~~~~~~~~(2.1)$$

\begin{thm}\label{thm1} Let $B,W$ be two independent standard Brownian motions. $W_0\equiv0.$
%$a,b,a',b',a_0,b_0,\beta$ are real numbers with the relationship $a<0<b,a\leq a'<b'\leq b,a< a_0<b_0< b.$
Under the probability $\bfP^{x},$ $B_{0}=x$ almost surely.
Define
$$\overline{X}_t:=-\ln\inf\limits_{x\in[a_0,b_0]} \bfP^{x}(\forall_{s\in [0,t]} \beta W_{s}+a\leq B_{s} \leq \beta W_{s}+b,~\beta W_{t}+a'\leq B_{t} \leq \beta W_{t}+b'|W);$$
$$\underline{X}_t:=-\ln ~\sup_{x\in\bfR}~\bfP^{x}(\forall_{s\in [0,t]} \beta W_{s}+a\leq B_{s} \leq \beta W_{s}+b|W).$$
Then there exists a function $\gamma: \bfR\rightarrow \bfR^{+}$ such that
 \beqlb\label{sect-1}
 \lim\limits_{t\rightarrow+\infty}\frac{\overline{X}_t}{t}=\lim\limits_{t\rightarrow+\infty}\frac{\underline{X}_t}{t}=\frac{\gamma(\beta)}{(b-a)^2},~~~{\rm a.s. ~~~and~~ in}~ L_p~(p\geq1),
 \eeqlb
 where $\gamma$ is a convex and even function. Moreover, $\gamma(0)=\frac{\pi^2}{2},\gamma(1)\leq 4\pi^2$ and for any $\beta\in\bfR,$ $\gamma(\beta)\geq \frac{\pi^2(1+\beta^2)}{2}.$ Hence
  $\gamma$ is strictly increasing on $[0,+\infty)$ and $\lim\limits_{|\beta|\rightarrow+\infty}\gamma(\beta)=+\infty.$
\end{thm}
\begin{rem}
In fact, (2.2) can be strengthened to (5.2) and (5.3). Moreover,  From the property of $\gamma,$ we can see even though the width of the random passageway is always``$b-a$'' at every moment from $0$ to $t,$ the first exit time will be shorter when the random passageway has more violent fluctuation of the center (i.e., when $|\beta|$ becomes bigger).
\end{rem}

In order to make writing more concise, we denote
\beqlb\label{sect-2}I^{W_s,\beta}_{~x,~y}:=[\beta W_s+x,\beta W_s+y],\eeqlb
where $x,y$ can be any constant or function.
\begin{cor} {\bf (Small deviation)}
If $\alpha\in(0,\frac{1}{2}),$   then we have, almost surely,
\beqlb\label{sect-3}
\liminf\limits_{t\rightarrow+\infty}\inf_{x\in[a_0t^\alpha,b_0t^\alpha]}\frac{\ln \bfP^x\big(\forall_{s\in [0,t]} B_s\in I^{W_s,\beta}_{at^{\alpha},bt^{\alpha}},B_{t}\in I^{W_t,~\beta}_{a't^{\alpha},b't^{\alpha}}|W\big)}{t^{1-2\alpha}}\geq-\frac{\gamma(\beta)}{(b-a)^2};~~~~~~~
\eeqlb
%~~~~~ (2.2). $$
\beqlb\label{sect-4}
\limsup\limits_{t\rightarrow+\infty}~\sup_{x\in\bfR}\frac{\ln \bfP^x(\forall_{s\in [0,t]} B_s\in[at^\alpha+\beta W_s, bt^\alpha+\beta W_s]|W)}{t^{1-2\alpha}}\leq-\frac{\gamma(\beta)}{(b-a)^2}.~~~~~\eeqlb
\end{cor}
\begin{rem} Obviously, the ``$\liminf, ~\geq $'' in (2.4) and ``$\limsup, ~\leq $'' in (2.5) can be replaced by ``$\lim,~=$''. The same replacement can also be done in (2.6) and (2.7).
\end{rem}
\begin{cor}
Let $f(s)$ and $g(s)$ be two continue functions from $[0,1]$ to $\bfR$ such that $$\forall s\in[0,1],~f(s)<g(s),~~f(0)<a_0\leq b_0<g(0),~~  f(1)\leq a'<b'\leq g(1).$$
We have, almost surely,
%$\inf_{x\in[0,1]}g(s)>\sup_{x\in[0,1]}$
\beqlb\label{sect-4}
 &&\liminf\limits_{t\rightarrow+\infty}\inf_{x\in[a_0,b_0]}\frac{\ln\bfP^{x}\Big(\forall_{s\in [0,t]}  B_{s}\in I^{~W_s,~\beta}_{f(\frac{s}{t}),g(\frac{s}{t})},B_{t}\in I^{W_t,\beta}_{~a',b'}|W\Big)}{t}\geq C_{f,g}\gamma(\beta),
\eeqlb\beqlb\label{sect-5}
\limsup\limits_{t\rightarrow+\infty}\sup_{x\in\bfR}\frac{\ln \bfP^{x}\big(\forall_{s\in [0,t]} \beta W_{s}+f(\frac{s}{t})\leq B_{s} \leq \beta W_{s}+g(\frac{s}{t})|W\big)}{t}\leq C_{f,g}\gamma(\beta),~~
\eeqlb
where $C_{f,g}:=-\int_{0}^{1}\big(g(s)-f(s)\big)^{-2}ds.$ We should notice that $C_{f,g}\in(0,+\infty)$ because of the assumption of $f(s)$ and $g(s)$.
%{\rm(}We should notice that $C_{f,g}\in(0,+\infty)$ because of the assumption of $f(s)$ and $g(s)$.{\rm)}
\end{cor}
%\begin{rem}
%In (3.2), it is regard to note that the `` living space " of the Brownian motion $B$ can be arbitrarily small only if $t$ is small enough, but for every $t$,
%$\bfP^0(\forall_{s\leq t} B_s\in[at^\alpha+\beta W_s, bt^\alpha+\beta W_s]|W)>0$ with probability 1.
%\end{rem}

\section{The moment estimation for $\overline{X}_t $}
The main tool we use to prove Theorem 2.1 is the Kingman's subadditive ergodic theorem. For preparation, we first give a important estimation for $\overline{X}_t$ which
has been defined in Theorem 2.1.

Since $\gamma$ is an even function and the distribution of the first exit time is well-known when $\beta=0$,
we will always assume $\beta>0$ in the rest of the paper.
\begin{thm}  For any $t>0, p>0, q>1,$ we have
\beqlb\label{sect3-1}\lim\limits_{n\rightarrow+\infty}n^p~ \bfP\big(\overline{X}_t \geq  (\ln n)^{q}\big)=0.~~\eeqlb
Thus for any $j\in\bfN,$ we have \beqlb\label{sect3-2} \bfE(\overline{X}^j_t)<+\infty.~~\eeqlb
\end{thm}
%\begin{rem} We suspect that there exists a constant $c_0\in\big(0,\frac{\pi^2}{2(b-a)^2}\big),$ such that for any $c<c_0, t>0,$ there has $\bfE(e^{c\overline{X}_t})<+\infty.$
%\end{rem}
{\bf Proof of Theorem 3.1}
Notice that $\overline{X}_t$ is related to $a,b,a_0,b_0,a',b'.$ Obviously, we only need to show that (3.1) holds when $a_0\leq a'< b'\leq b.$
Recalling the basic relationship (2.1), we will first prove (3.1) under the situation of
\beqlb\label{sect3-3} a_0\leq a'< b'\leq b~~~ \text{and}~~~ \min\{a_0-a,b-b_0\}>\max\{a'-a_0,b_0-b'\}.\eeqlb
Under this situation, we can choose $a'',b'',\delta$ such that $a'<a''<b''<b'$ and
\beqlb\label{sect3-4}0<2\beta\delta < \min\big\{b''-a'',~\min\{a_0-a,b-b_0\}-\max\{a''-a_0,b_0-b''\} \big\}.\eeqlb
We define a Markov time sequence $\{\tau_{_{n,\delta}}\}_{n\in\bfN}$ such that
$$\tau_{_{0,\delta}}:=0,~~~ \tau_{_{n+1,\delta}}:=\inf\{s>\tau_{_{n,\delta}}: |W_{s}-W_{\tau_{_{n,\delta}}}|=\delta\}.~~~~n=0,1,2,\ldots$$
It is easy to see that $\{\tau_{_{n,\delta}}\}_{n\in\bfN}$ is an i.i.d. random walk and $\tau_{1,\delta}>0$ almost surely.
We divide time $[0,t]$ into $[0,\tau_{_{1,\delta}}], [\tau_{_{1,\delta}},\tau_{_{2,\delta}}],...[\tau_{_{N,\delta}}, t],$  where
$$N:= \sup\{n: \tau_{_{n,\delta}}< t\}.$$
Then by the Markov property, we have
\beqnn\overline{X}_t\leq \mathbf{1}_{\{N=0\}}Z_{0}(W)+\sum_{i=1}^{+\infty}\mathbf{1}_{\{N=i\}}\Big(\sum_{k=0}^{i-1}Y_k(W)+Z_{i}(W)\Big), \eeqnn
where
\beqnn &Y_k(W):=-\ln&\inf_{x\in[a_0,b_0]}\bfP\big(\forall_{\tau_{_{k,\delta}}\leq s\leq \tau_{_{k+1,\delta}}} B_{s}-\beta (W_s- W_{\tau_{_{k,\delta}}})\in[a,b],
\\&~~~~~~~~~&B_{\tau_{_{k+1,\delta}}}-\beta (W_{\tau_{_{k+1,\delta}}}-W_{\tau_{_{k,\delta}}})\in[a'',b'']~|W, B_{\tau_{_{k,\delta}}}=x\big),\eeqnn
$$Z_k(W):=-\ln\inf_{x\in[a'', b'']}\bfP(\forall_{\tau_{_{k,\delta}}\leq s\leq t}B_{s}-(\beta W_s-\beta W_{\tau_{_{k,\delta}}})\in[a',b']~|W,B_{\tau_{_{k,\delta}}}=x).$$
Let $$\rho_{_{k,\delta}}:=\tau_{_{k,\delta}}-\tau_{_{k-1,\delta}},~k=1,2,3,\ldots$$
By the definition of $\tau_{_{k,\delta}}$ we can get a further upper bound for $Y_k(W)$ and $Z_k(W),$ which is
 \beqnn Y_k&:=&-\ln\inf_{x\in[a_0,b_0]}\bfP^{x}(\forall_{0\leq s\leq \rho_{_{k+1,\delta}}}B_{s}\in[a+\beta\delta,b-\beta\delta],B_{\rho_{_{k+1,\delta}}}\in[a''+\beta\delta,b''-\beta\delta]|W)
 \\&\geq& Y_k(W),\eeqnn $$ Z_0:=-\ln\inf_{x\in[a'',b'']}\bfP^{x}(\forall_{0\leq s\leq t}B_{s}\in[a'+\beta\delta,b'-\beta\delta])\geq Z_k(W).$$
 Note that $Y_k$ is also depend on $W$ but $Z_0$ is a non-random constant. Hence we have
\beqlb\label{sect3-5}\overline{X}_t\leq \mathbf{1}_{\{N=0\}}Z_{0}+\sum_{i=1}^{+\infty}\mathbf{1}_{\{N=i\}}(\sum_{k=0}^{i-1}Y_k+Z_0)=Z_0+\sum_{i=1}^{+\infty}\mathbf{1}_{\{N=i\}}\sum_{k=0}^{i-1}Y_k. \eeqlb
Naturally, we need to estimate the upper bound of $Y_k.$
Define \beqlb\label{sect3-6} k(t):=\inf_{x\in[a_0,b_0]}\bfP^x(\forall_{0\leq s\leq t}B_s\in[a+\beta\delta,b-\beta\delta],B_{t}\in[a''+\beta\delta,b''-\beta\delta]),\eeqlb
$$\delta_1:= \max\{a''+\beta\delta-a_0,b_0-b''+\beta\delta\},~~ \delta_2:=\min\{a_0-a-\beta\delta,b-\beta\delta-b_0\}.$$
By basic calculation, we have
\beqnn
k(t)&\geq& \inf_{x\in[a_0,b_0]}\Big\{\bfP^x(B_{t}\in[a''+\beta\delta,b''-\beta\delta])-\big[1-\bfP^x\big(\forall_{0\leq s\leq t}B_s\in\left[a-\beta\delta,b+\beta\delta\right]\big)\big]\Big\}
\\&\geq& \bfP^{0}(B_{t}\in[\delta_1,\delta_1+b''-a''-2\beta\delta])-\bfP^0(\sup\limits_{s\in[0,t]}|B_s|>\delta_2).
\eeqnn
By (\ref{sect3-4}), we know $\delta_2> \delta_1.$  Therefore, we can choose an $\epsilon>0$ small enough such that
$\frac{(\delta_1+\epsilon)^2}{2}\leq \frac{\delta_2^2}{2+\epsilon}$ and $\delta_1+\epsilon\leq \delta_1+b''-a''-2\beta\delta.$
That means
$$\bfP^{0}(B_{t}\in[\delta_1,\delta_1+b''-a''-2\beta\delta])\geq \frac{\epsilon}{\sqrt{2\pi t}}\exp\Big\{-\frac{(\delta_1+\epsilon)^2}{2t}\Big\}.$$
Recalling the Csorgo and Revesz estimation \cite[Lemma1]{CR1979}, we know there exists a constant $C>0$ such that
$$\bfP^0(\sup\limits_{s\in[0,t]}|B_s|>\delta_2)\leq C\exp\Big\{-\frac{\delta_2^2}{(2+\epsilon)t}\Big\}.$$
Then there exists a $D>0$ such that for any $t\leq D,$
\beqlb\label{sect3-7}
k(t)\geq \frac{\epsilon}{\sqrt{2\pi t}}\exp\Big\{-\frac{(\delta_1+\epsilon)^2}{2t}\Big\}-C\exp\Big\{-\frac{\delta_2^2}{(2+\epsilon)t}\Big\}\geq \exp\Big\{-\frac{(\delta_1+\epsilon)^2}{2t}\Big\}.~~~~\eeqlb
When $t>D,$ by (\ref{sect3-4}), we can choose a $\delta_3>0$ such that
$a'+\beta\delta+\delta_3< b'-\beta\delta-\delta_3.$
So for $t>D$, we have
\beqlb\label{sect3-8}
k(t)&\geq&\inf_{x\in[a_0,b_0]}\bfP^x(\forall_{0\leq s\leq D}B_s\in[a+\beta\delta,b-\beta\delta],B_{D}\in[a'+\beta\delta+\delta_3,b'-\beta\delta-\delta_3])\nonumber
\\&\times&\inf_{x\in[a'+\beta\delta+\delta_3,b'-\beta\delta-\delta_3]}\bfP^x(\forall_{0\leq s\leq t-D}B_s\in[a'+\beta\delta,b'-\beta\delta]).
\eeqlb
Notice that
$\inf_{x\in[a_0,b_0]}\bfP^x(\forall_{0\leq s\leq D}B_s\in[a+\beta\delta,b-\beta\delta],B_{D}\in[a'+\beta\delta+\delta_3,b'-\beta\delta-\delta_3])$
is a positive constant. Moreover,
\beqlb\label{sect3-9}
&&\inf_{x\in[a'+\beta\delta+\delta_3,b'-\beta\delta-\delta_3]}
\bfP^x\big(\forall_{0\leq s\leq t-D}B_s\in[a'+\beta\delta,b'-\beta\delta]\big)\nonumber
\\&&\geq\bfP^0(\sup\limits_{s\in[0,t-D]}|B_s|\leq \delta_3)\geq\frac{8}{3\pi}\exp\Big\{-\frac{\pi^{2}(t-D)}{8\delta_3^2}\Big\}.
\eeqlb
Combining with (\ref{sect3-7})-(3.9), we conclude that there exist $C_1,C_2>0$ such that
\beqlb\label{sect3-10}-\ln k(t)\leq {C_1}t^{-1}\mathbf{1}_{\{t\leq D\}}+C_2t\mathbf{1}_{\{t> D\}}.\eeqlb
It implies that for any $k\in\{1,2,\ldots,N\},$ we have
\beqlb\label{sect3-11}Y_{k-1}\leq {C_1}\rho^{-1}_{_{k,\delta}}\mathbf{1}_{\{\rho_{_{k,\delta}}\leq D\}}+C_2\rho_{_{k,\delta}}\mathbf{1}_{\{\rho_{_{k,\delta}}> D\}}.\eeqlb
Choosing $q'',q'$ such that $1<q''<q'<q.$ When $n$ is large enough, by (\ref{sect3-5}) and (3.11)  we have
\beqnn\bfP(\overline{X}_t\geq (\ln n)^{q})&\leq& \sum_{i=1}^{+\infty}\bfE(\mathbf{1}_{\{\sum_{k=0}^{i-1}Y_i+Z_0\geq (\ln n)^{q}\}}\mathbf{1}_{\{N=i\}})
+ \bfE(\mathbf{1}_{\{Z_0\geq (\ln n)^{q}\}}\mathbf{1}_{\{N=0\}})
\\&\leq& \sum_{i=1}^{+\infty}\bfE(\mathbf{1}_{\{\sum_{k=0}^{i-1}Y_k\geq 2(\ln n)^{q'}\}}\mathbf{1}_{\{N=i\}})
%\\&\leq& \sum_{i=1}^{+\infty}\bfP\Big(\sum_{k=1}^{i}\big(\frac{C_1}{\rho_{k,\delta}}+C_2\rho_{k,\delta}\big)\geq 2(\ln n)^{q'},N=i\Big)
\\&\leq& \sum_{i=1}^{+\infty}\Big[\bfE\big(\mathbf{1}_{\{\sum_{k=1}^{i}\frac{C_1}{\rho_{_{k,\delta}}}\geq (\ln n)^{q'},N=i\}}\big)+\bfE\big(\mathbf{1}_{\{\sum_{k=1}^{i}C_2\rho_{_{k,\delta}}\geq (\ln n)^{q'},N=i\}}\big)\Big].
\eeqnn
Notice that when $(\ln n)^{q'}\geq C_2t,$ for any $i,$ we have
$$\bfP\big(\sum_{k=1}^{i}C_2\rho_{_{k,\delta}}\geq (\ln n)^{q'},N=i\big)\leq\bfP\big(\tau_{_{i,\delta}}\geq\frac{(\ln n)^{q'}}{C_2}, \tau_{_{i,\delta}}< t\big)=0.$$
Hence when $n$ is large enough, denote $\varsigma_n:=(\ln n)^{q''},$ it is true that
\beqlb\label{sect3-12}\bfP(\overline{X}_t\geq (\ln n)^{q})&\leq& \sum_{i=1}^{+\infty}\bfP\Big(\sum_{k=1}^{i}\frac{C_1}{\rho_{_{k,\delta}}}\geq (\ln n)^{q'},N=i \Big)\nonumber
\\&\leq& \sum_{i=1}^{\lfloor\varsigma_n\rfloor-1}\bfP\Big(\sum_{k=1}^{i}\frac{C_1}{\rho_{_{k,\delta}}}\geq (\ln n)^{q'} \Big)+\sum_{i=\lfloor\varsigma_n\rfloor}^{+\infty}\bfP(N=i)\nonumber
\\&\leq& \varsigma_n\bfP\Big(\sum_{k=1}^{\lfloor\varsigma_n\rfloor}\frac{1}{\rho_{_{k,\delta}}}\geq \frac{(\ln n)^{q'}}{C_1}\Big)+
\bfP(\tau_{_{\lfloor\varsigma_{_n}\rfloor,\delta}}\leq t),
\eeqlb
where $\lfloor x\rfloor$ is the integer part of $x$. Let $p(t)$ be the probability density function of Markov time $\tau_{_{1,\delta}}.$ Then the expression of $p(t)$ is
$$p(t)=\frac{2\delta}{\sqrt{2\pi t^3}}\sum\limits_{n=-\infty}^{+\infty}(4n+1)e^{-\frac{(4n+1)^2\delta^2}{2t}}.$$
Obviously, $p(t)\leq  \frac{2\delta}{\sqrt{2\pi t^3}}e^{-\frac{\delta^2}{2t}}$ when $t$ is small enough.
 Then it is easy to see for any $c_3>0,$ we have $\bfE(\exp\{\frac{c_3}{\tau_{_{1,\delta}}}\})<+\infty.$ Moreover, according to \cite[Page 30]{IM1974}, there exists a positive constant $c_4>0$ such that for any $c_5\in[0,c_4]$, $\bfE(e^{c_{_5}\tau_{_{1,\delta}}})<+\infty$. Then we can apply the Cram\'{e}r theorem \cite[Page 27]{DZ1998} to i.i.d. random walk $\{\frac{1}{\rho_{_{1,\delta}}}+\frac{1}{\rho_{_{2,\delta}}}+\cdots+\frac{1}{\rho_{_{i,\delta}}}\}_{i\in\bfN}$ and $\{\tau_{_{i,\delta}}\}_{i\in\bfN}$ respectively. Notice that
$\bfE(\tau_{_{1,\delta}}),\bfE(\frac{1}{\tau_{_{1,\delta}}})\in(0,+\infty),q''<q'$ and time $t$ is not depend on $n,$ so when $n$ is large enough there exist $c_1,c_2>0$ such that
$$\bfP\Big(\sum_{k=1}^{\lfloor\varsigma_n\rfloor}\frac{1}{\rho_{_{k,\delta}}}\geq \frac{(\ln n)^{q'}}{C_1}\Big)
\leq\bfP\Big(\frac{\sum_{k=1}^{\lfloor(\ln n)^{q''}\rfloor}\rho^{-1}_{_{k,\delta}}}{\lfloor(\ln n)^{q''}\rfloor}\geq 2\bfE\big(\frac{1}{\tau_{_{1,\delta}}}\big)\Big) \leq e^{-c_1\lfloor(\ln n)^{q''}\rfloor}$$
and
$$\bfP(\tau_{_{\lfloor\varsigma_{_n}\rfloor,\delta}}\leq t)\leq
\bfP\Big(\frac{\tau_{\lfloor(\ln n)^{q''}\rfloor,\delta}}{\lfloor(\ln n)^{q''}\rfloor}\leq \frac{\bfE(\tau_{_{1,\delta}})}{2}\Big)\leq e^{-c_2\lfloor(\ln n)^{q''}\rfloor}.$$
Combining with (3.12), we get, as $n\rightarrow +\infty,$
$$n^p\bfP(\overline{X}_t\geq (\ln n)^{(1+q)})\leq n^p(\ln n)^{q''}e^{-c_1\lfloor(\ln n)^{q''}\rfloor}+n^pe^{-c_2\lfloor(\ln n)^{q''}\rfloor}\rightarrow 0.$$
So we obtain (3.1) under situation (3.3).

Next, if (3.3) is not hold, that is to say, $\min\{a_0-a,b-b_0\}\leq\max\{a'-a_0,b_0-b'\}.$
Let $\lceil x\rceil:=\inf\{j\in\bfN: j\geq x\},$ define $m:=\Big\lceil\frac{3\max\{a'-a_0,b_0-b'\}}{\min\{a_0-a,b-b_0\}}\Big\rceil,~~ \delta':=\frac{\min\{a_0-a,b-b_0\}}{2}, $ $a_i:=\min\{a_0+i\delta',a'\},~b_i:=\max\{b_0-i\delta',b'\},$ and
\beqnn\overline{X}_{t,i}:=-\ln\inf\limits_{x\in[a_i,b_i]} \bfP(\forall_{s\in [\frac{it}{m},\frac{(i+1)t}{m}]} B_{s}-\beta (W_{s}-W_{\frac{it}{m}})\in [a,b],~~~~~~~~
\\~~~~~~~~~B_{\frac{(i+1)t}{m}}-\beta (W_{\frac{(i+1)t}{m}}-W_{\frac{it}{m}})\in [a_{i+1},b_{i+1}]|W, B_{\frac{it}{m}}=x).\eeqnn
By the Markov property we can see that $\overline{X}_t \leq \sum_{i=0}^{m-1}\overline{X}_{t,i}.$
%$$\overline{X}_{t,i}:=-\ln\inf\limits_{x\in[a_i,b_i]} \bfP\Big(\substack{\forall_{s\in [\frac{it}{m},\frac{(i+1)t}{m}]} B_{s}-\beta (W_{s}-W_{\frac{it}{m}})\in [a,b] \\ B_{\frac{(i+1)t}{m}}-\beta (W_{\frac{(i+1)t}{m}}-W_{\frac{it}{m}})\in [a_{i+1},b_{i+1}]}\Big|W, B_{\frac{it}{m}}=x\Big).$$

Notice that $\overline{X}_{t,i}$ has the same law as $\overline{X}'_{t,i},$ where
\beqnn\overline{X}'_{t,i}:=-\ln\inf\limits_{x\in[a_i,b_i]} \bfP^x(\forall_{s\in [0,\frac{t}{m}]} B_{s}-\beta W_{s}\in [a,b],
B_{\frac{t}{m}}-\beta W_{\frac{t}{m}}\in [a_{i+1},b_{i+1}]|W).\eeqnn
As $q'<q,$ we have $(\ln n)^q\geq m(\ln n)^{q'}$ when $n$ is large enough. Consequently,
\beqlb\label{sect3-13}n^p\bfP(\overline{X}_t\geq (\ln n)^{q})\leq n^p\sum_{i=0}^{m-1}\bfP(\overline{X}_{t,i}\geq (\ln n)^{q'})=\sum_{i=0}^{m-1}n^p\bfP(\overline{X}'_{t,i}\geq (\ln n)^{q'}).~~~~~~\eeqlb
%\\&\leq& (\ln n)^{q''}e^{-c_1\lfloor(\ln n)^{q''}\rfloor}+e^{-c_2\lfloor(\ln n)^{q''}\rfloor}.
 Note that for any $i\in[0,m-1]\cap\bfN,$ $a,b,a_i,b_i,a_{i+1},b_{i+1}$ satisfy the relationship (3.3), so we have $$\lim\limits_{n\rightarrow 0}n^p\bfP(\overline{X}'_{t,i}\geq (\ln n)^{q'})=0,~~~i\in [0,m-1]\cap\bfN.$$
Combining with (3.13), we complete the proof of (3.1).

Note that for any $j\in\bfN,$
$$\bfE(\overline{X}^j_t)\leq \sum_{n=0}^{+\infty}(n+1)^j\bfP(\overline{X}^j_t\in [n,n+1])\leq \sum_{n=0}^{+\infty}(n+1)^j\bfP(\overline{X}^j_t\geq n).$$
Moreover, for large enough $n,$ according to (3.1) we have $\bfP(\overline{X}^j_t\geq n)\leq e^{-\sqrt{n}},$ which implies that (3.2) holds.

\qed

\section{Proof of the main result}
In this section, we will show how to use the Kingman's subadditive ergodic theorem \cite[Theorem 9.14]{K1997} to prove Theorem 2.1. To simplify the statement,we first introduce some notations.
Let $t_2>t_1\geq 0,$  analogous to the definition of (2.3), we define
$$I^{W_s,\beta}_{a,b,t_1}:=[a+\beta W_s-\beta W_{t_1},b+\beta W_s-\beta W_{t_1}].$$
Denote ~$r_{t_1,t_2}(a,b,a',b',x,\beta):=\bfP(\forall_{t_1\leq s\leq t_2}B_{s}\in I^{W_s,\beta}_{a,b,t_1},B_{t_2}\in I^{W_{t_2},\beta}_{a',b',t_1}|W,B_{t_1}=x),$
$$p_{t_1,t_2}(a,b,a',b',\beta):=\inf_{x\in[a',b']}r_{t_1,t_2}(a,b,a',b',x,\beta),~~$$
$$q_{t_1,t_2}(a,b,a',b',\beta):=-\ln p_{t_1,t_2}(a,b,a',b',\beta).$$
Without causing confusion, sometimes $p_{t_1,t_2}(a,b,a',b',\beta)$ and $q_{t_1,t_2}(a,b,a',b',\beta)$ are abbreviated as $p_{t_1,t_2}$ and $q_{t_1,t_2}$ respectively in the rest part of the paper.
The following lemma is essential for Theorem 2.1.
\begin{lem}
 Under the situation $a<a'<b'<b$ and the relationship (2.1), there exists a non-negative function of two variables $\gamma: (0,+\infty)\times\bfR\rightarrow [0,+\infty)$ such that
 \beqlb\label{sect3-9}
\lim\limits_{n\rightarrow+\infty} \frac{q_{0,n}(a,b,a',b',\beta)}{n}=\gamma(b-a,\beta),    ~~~    {\rm a.s.~~and~~in}~~L^1.
 \eeqlb
% for every $x\in(a,b)$
%$$-\lim\limits_{n\rightarrow+\infty}\frac{\ln r_{0,n}(a,b,a',b',x,\beta)}{n}=\gamma(a,b,\beta)  \   ~~~     a.s.$$
\end{lem}

\noindent{\bf Proof of Lemma 4.1}.  We divide the proof into two steps.

\smallskip

\noindent\emph{Step 1. Showing that $\{\frac{q_{0,n}(a,b,a',b',\beta)}{n}\}$ has an almost surely and $L^1$ degenerate limit.}

By the Markov property, we know
\beqnn p_{0,n}&\geq& p_{0,m}\inf_{x\in I^{W_m,\beta}_{~a',~b'}}\bfP(\forall_{m\leq s\leq n}B_{s}\in I^{W_s,\beta}_{~a,~b},B_{n}\in I^{W_n,\beta}_{~a',~b'}|W,B_{m}=x)
\\&=& p_{0,m}p_{m,n},  ~~~~~~~0\leq m< n. \eeqnn Hence we have $q_{0,n}\leq q_{0,m}+q_{m,n}.$
 This is the subadditivity condition \cite[(9.9)]{K1997} of the Kingman's subadditive ergodic theorem.

If we denote $W^i(s):=W_{i+s}-W_i, s\in[0,1],$  it is easy to see the sequence $\{W^i\}_{i\in \bfN}$ is i.i.d. and the randomness of $q_{m,n}$ is only depend on
$\{W_{m+s}-W_m, s\in[0,n-m]\}.$
From these facts and the stationary independent increments property of Brownian motion, we know that for any fixed $k$, the random sequence $q_{0,k},q_{k,2k},...,q_{nk,(n+1)k},...$ is i.i.d., and for every $l\in\bfN$, random sequence ${q_{l,l+1},q_{l,l+2},...,q_{l,l+n},...}$ has the same distribution as ${q_{0,1},q_{0,2},...,q_{0,n},...}$. These mean that $\{q_{m,n}\}_{1\leq m\leq n}$ fulfills the conditions \cite[(9.7)]{K1997} and \cite[(9.8)]{K1997} respectively. According to Theorem 3.1, we know $\bfE(q_{0,1})<+\infty,$ which is the
integrability condition of the Kingman's subadditive ergodic theorem. And obviously, for
each $n$, $\frac{\bfE(q_{0,1})}{n}\geq 0>-\infty.$  So far we have verified all conditions of the Kingman's subadditive ergodic theorem.

Besides, for every $k$, the sequence ${q_{0,k},q_{k,2k},...,q_{nk,(n+1)k},...}$ is ergodic since it is i.i.d., and thus we can conclude that $\frac{q_{0,n}}{n}$ converges to a constant almost surely and in $L_1$. Here we denote the limit by $\gamma(a,b,a',b',\beta)$. Consequently, we have
$$\lim\limits_{n\rightarrow +\infty}\frac{q_{0,n}(a,b,a',b',\beta)}{n}=\gamma(a,b,a',b',\beta),~~~    {\rm a.s.~~and~~in}~~L^1.$$
    %By the way, according to the Kingman's subadditive theorem, the convergence also makes means in the sense of $L^1$, but that is not involved in the research of this paper, so we omit it.

    \smallskip

\noindent\emph{Step 2.~~Let $a<a''<b''<b,$ showing that $\gamma(a,b,a',b',\beta)=\gamma(a,b,a'',b'',\beta).$}

Without loss of generality, we assume $a'\leq a''\leq b''\leq b'.$
Obviously,
%To prove this, we only need to show that for $\forall a'\leq a''\leq b''\leq b'$,we have $$\gamma(a,b,a',b',\beta)=\gamma(a,b,a'',b'',\beta)$$
%Now we do it.
\beqlb\label{sect4-2} q_{0,n}(a,b,a'',b'',\beta)&\leq& -\ln \inf\limits_{x\in[a',b']}r_{0,n}(a,b,a'',b'',x,\beta)\nonumber
\\&\leq& -\ln \inf\limits_{x\in[a',b']}r_{0,1}(a,b,a'',b'',x,\beta)+q_{1,n}(a,b,a'',b'',\beta). \eeqlb
By step {\emph 1} and the stationary increments property of Brownian motion, we know $\frac{q_{1,n}(a,b,a'',b'',\beta)}{n-1}\rightarrow \gamma(a,b,a'',b'',\beta)$ in probability. Moreover, applying the Kingman's subadditive ergodic theorem again, we can see $\frac{q_{1,n}(a, b, a'',b'',\beta)}{n-1}$ converges to a constant almost surely. Hence we have  $$\frac{q_{1,n}(a,b,a'',b'',\beta)}{n-1}\rightarrow \gamma(a,b,a'',b'',\beta),~~~~~{\rm a.s.}.$$
According to Theorem 3.1 and the Borel-Cantelli 0-1 law, it is easy to see that for any function $\alpha: \bfN\rightarrow \bfN$ and $k\in\bfN,$  \beqlb\label{sect4-3} \lim\limits_{n\rightarrow+\infty}\frac{-\ln\inf_{x\in[a',b']}r_{\alpha(n),\alpha(n)+k}(a,b,a',b',x,\beta)}{n}=0,~~~~~{\rm a.s.}.\eeqlb Hence from (\ref{sect4-2}) and (\ref{sect4-3}) we can obtain
\beqlb\label{sect4-4} \lim\limits_{n\rightarrow+\infty}\frac{-\ln \inf_{x\in[a',b']}r_{0,n,}(a,b,a'',b'',x,\beta)}{n}=\gamma(a,b,a'',b'',\beta),~~~~ {\rm a.s.} .\eeqlb
On the other hand, we have
\beqnn q_{0,n}(a,b,a',b',\beta)&\leq& -\ln \inf\limits_{x\in[a',b']}r_{0,n}(a,b,a'',b'',x,\beta)
\\&\leq& q_{0,n-1}(a,b,a',b',\beta)-\ln \inf\limits_{x\in[a',b']}r_{n-1,n}(a,b,a'',b'',x,\beta). \eeqnn
%$$ \lim\limits_{n\rightarrow+\infty}-\ln \inf\limits_{x\in[a',b']}\frac{r_{0,n,\delta}(a,b,a'',b'',x,\beta)}{n}=\gamma_\delta(a,b,a',b',\beta)~~~~~~~~ a.s.$$
Analogous to the above discussion, we can also obtain
$$ \lim\limits_{n\rightarrow+\infty}\frac{-\ln \inf_{x\in[a',b']}r_{0,n}(a,b,a'',b'',x,\beta)}{n}=\gamma(a,b,a',b',\beta),~~~~~~~~ {\rm a.s.}.$$
Combining with (\ref{sect4-4}), we have $\gamma(a,b,a',b',\beta)=\gamma(a,b,a'',b'',\beta).$

%That implies for any fixed pair $(a',b')$ satisfied that $(a'-a-\beta\delta)\wedge(b-\beta\delta-b')>\beta\delta,$
%if $$a''\in [a',a'-\beta\delta+(a'-a-\beta\delta)\wedge(b-\beta\delta-b')), b''\in [b'+\beta\delta-(a'-a-\beta\delta)\wedge(b-\beta\delta-b'),b'),$$
%(Here we should note it always can find a small enough $\delta$ such that $a'<a'',b''<b'$)
%then we always have $\gamma_\delta(a,b,a',b',\beta)=\gamma_\delta(a,b,a'',b'',\beta).$

%Take same steps, we can know if $(a,b,a'',b'',a''',b''',\delta)$ satisfy condition 1. we can obtain $\gamma_\delta(a,b,a'',b'',\beta)=\gamma_\delta(a,b,a''',b''',\beta).$ It means that for any $a'\leq a''\leq b''\leq b'$, $\gamma_\delta(a,b,a',b',\beta)=\gamma_\delta(a,b,a'',b'',\beta).$

If $a',b',a'',b''$ can not satisfy the relationship  `` $a'\leq a''\leq b''\leq b'$ ", then without loss of generality, we assume $a<a'\leq a''\leq b'\leq b''<b.$ From the above conclusion, we have
$$\gamma(a,b,a',b',\beta)=\gamma(a,b,a'',b',\beta)=\gamma(a,b,a'',b'',\beta).$$
So it is reasonable to write $\gamma(a,b,a',b',\beta)$ as $\gamma(a,b,\beta).$ Moreover, by the basic property of Brownian motion, it is easy to see for any $c\in\bfR,$ we have $\gamma(a+c,b+c,\beta)=\gamma(a,b,\beta).$ Hence we can further denote $\gamma(a,b,\beta)$ by $\gamma(b-a,\beta).$
This is the end of the proof of Lemma 4.1.

\qed

More information of $\gamma(c,\beta)$ has been listed in the following proposition, which is also a important preparation for the proof of Theorem 2.1.

\begin{prop} The function $\gamma(c,\beta)$ has been introduced in Lemma 4.1.
\item[(1).]  For each fixed $\beta,$ the function $c\mapsto\gamma(c,\beta)$ is convex on $(0,+\infty)$.
\item[(2).]  For each fixed $c>0,$ the function $\beta\mapsto\gamma(c,\beta)$ is even and convex.
\end{prop}

\noindent{\bf Proof of proposition 4.1}.
%The main idea of this proof is utilizing the log-concave property of Gaussian measure. Now we show the detail.

In this proof, we set $a,b,a',b'$ satisfy the relationship 2.1 and $a<a'\leq0\leq b'<\min \{b_1, b_2, b\}.$ Denote $d_{n,T}(x)$ is the joint density function of $(B_{_{T/n}},B_{_{2T/n}}...,B_{_{nT/n}})$ from $\bfR^{n}$ to $\bfR^{+}.$ By basic calculation, we know for any $n,T,$ function $d_{n,T}$ is a log-concave function, that is to say, for any $\lambda\in[0,1]$ and $x,y\in \bfR^{n},$ it has the relationship $d_{n,T}(\lambda x+(1-\lambda) y)\geq d_{n,T}(x)^\lambda d_{n,T}(y)^{(1-\lambda)}.$

Let $g_1(s),g_2(s),h_1(s),h_2(s)$ be real functions defined on $[0,1]$ such that for every $s\in[0,1]$, $g_1(s)\leq h_1(s),~g_2(s)\leq h_2(s).$
If we denote the $k$-th coordinate of $x,y\in \bfR^{n}$ by $x_k,y_k$, it is obvious that
\beqnn&& 1_{\forall_{k\leq n} \lambda x_k+(1-\lambda)y_k\in[\lambda g_{1}(k/n)+(1-\lambda) g_{2}(k/n),\lambda h_{1}(k/n)+(1-\lambda) h_{2}(k/n)]}
\\&\geq&\big(1_{\forall_{k\leq n}  x_k\in[g_{1}(k/n),h_{1}(k/n)]}\big)^{\lambda}\big(1_{\forall_{k\leq n}  y_k\in[g_{2}(k/n),h_{2}(k/n)]}\big)^{1-\lambda}.\eeqnn
Denote $H^{\lambda}_{n,T}(x):=d_{n,T}(x)1_{\forall_{k\leq n} x_k\in[\lambda g_{1}(k/n)+(1-\lambda) g_{2}(k/n),\lambda h_{1}(k/n)+(1-\lambda) h_{2}(k/n)]},$
then we have
\beqlb\label{sect4-5} H^{\lambda}_{n,T}(\lambda x+(1-\lambda) y)\geq (H^{1}_{n,T}(x))^{\lambda} (H^{0}_{n,T}(y))^{1-\lambda}.\eeqlb
Moveover, by Theorem 3.1, we know for each $m\in\bfN,$
\beqlb\label{sect4-6}q_{0,m}(a,b,a',b',\beta)<+\infty,~~~~~{\rm a.s.}.\eeqlb
 (\ref{sect4-5}) and (\ref{sect4-6}) are the two conditions of the Prekopa-Leindler inequality \cite[Theorem 7.1]{G2002}.
 According to the Prekopa-Leindler inequality, we have
 \beqlb\label{sect4-7}\int_{\mathbb{R}^n}H^{\lambda}_{n,m}(x)dx\geq\Big(\int_{\mathbb{R}^n}H^1_{n,m}(x)dx\Big)^{1-\lambda}\Big(\int_{\mathbb{R}^n}H^0_{n,m}(x)dx\Big)^\lambda.\eeqlb
 If we set
 $$g_{i}(s):=\begin{cases}\beta W_{sm}+a,&{s\in[0,1)}
 \\\beta W_{sm}+a',&{s=1},\end{cases}{\rm and}~~h_{i}(s):=\begin{cases}\beta W_{sm}+b_i,&{s\in[0,1)}
 \\\beta W_{sm}+b',&{s=1},\end{cases} ~~i=1,2.$$
 then (\ref{sect4-7}) means that except the zero measure set $\{\omega: q_{0,m}(a,b,a',b',\beta)=+\infty\},$ we always have
 \beqnn
 &&\bfP\Big(\forall_{1\leq k\leq n,k\in\bfN}B_{km/n}\in I^{W_{km/n},\beta}_{a,\lambda b_1+(1-\lambda)b_2},B_m\in I^{W_m,\beta}_{a',b'}|W,B_0=0\Big)
 \\&\geq&\Big[\bfP\Big(\forall_{1\leq k\leq n,k\in\bfN}B_{km/n}\in I^{W_{km/n},\beta}_{a,b_1},B_m\in I^{W_m,\beta}_{a',b'}|W,B_0=0\Big)\Big]^\lambda~~~
 \\&~~&\times\Big[\bfP\Big(\forall_{1\leq k\leq n,k\in\bfN}B_{km/n}\in I^{W_{km/n},\beta}_{a,b_2},B_m\in I^{W_m,\beta}_{a',b'}|W,B_0=0\Big)\Big]^{1-\lambda}.
 \eeqnn
Let $n\rightarrow +\infty,$ we deduce that for each $m\in \mathbb{N}$, almost surely we have
 $$r_{0,m}(a,\lambda b_1+(1-\lambda)b_2,a',b',0,\beta)\geq r^\lambda_{0,m}(a,b_1,a',b',0,\beta)r^{1-\lambda}_{0,m}(a,b_2,a',b',0,\beta).$$
For any $c\in[a',b'],$ by the same way, we can prove
 \beqnn
 &&r_{0,m}(a-c,\lambda b_1+(1-\lambda)b_2-c,a'-c,b'-c,0,\beta)
 \\&\geq& r^\lambda_{0,m}(a-c,b_1-c,a'-c,b'-c,0,\beta)\times r^{1-\lambda}_{0,m}(a-c,b_2-c,a'-c,b'-c,0,\beta),
 \eeqnn
which means that
 \beqnn
 &&\inf_{x\in[a',b']}r_{0,m}(a,\lambda b_1+(1-\lambda)b_2,a',b',x,\beta)
 \\&\geq&\inf_{x\in[a',b']}\big(r^\lambda_{0,m}(a,b_1,a',b',x,\beta)\times r^{1-\lambda}_{0,m}(a,b_2,a',b',x,\beta)\big)
 \\&\geq&\inf_{x\in[a',b']}r^\lambda_{0,m}(a,b_1,a',b',x,\beta)\times\inf_{y\in[a',b']} r^{1-\lambda}_{0,m}(a,b_2,a',b',y,\beta),~~~~{\rm a.s.}.
 \eeqnn
That is to say
 \beqlb\label{sect-13} &&q_{0,m}(a,\lambda b_1+(1-\lambda)b_2,a',b',\beta)\nonumber
 \\&\leq&\lambda q_{0,m}(a,b_1,a',b',\beta)+(1-\lambda)q_{0,m}(a,b_2,a',b',\beta).~~~~{\rm a.s.}.\eeqlb
Therefore, by Lemma 4.1 we can see the function $c\mapsto\gamma(c,\beta)$ is convex on $(0,+\infty)$.

Now it is time to show Proposition 4.1 (2). Obviously, the function $\beta\mapsto\gamma(c,\beta)$ is even since the standard Brownian motion $W$ is symmetric. Hence we only need to consider $\beta\in[0,+\infty).$ For any $\beta_1,\beta_2\geq 0$, if we set
 $$g_{i}(s):=\begin{cases} \beta_i W_{sm}+a, &{s\in[0,1)}
 \\\beta_i W_{sm}+a', &{s=1},\end{cases} {\rm and}~~h_{i}(s):=\begin{cases} \beta_i W_{sm}+b, &{s\in[0,1)}
 \\\beta_i W_{sm}+b', &{s=1},\end{cases} ~~i=1,2,$$
then we can obtain Proposition 4.1 (2) by repeating the step (4.5)-(4.8) similarly.

\qed

Now we will prove Theorem 2.1.

\noindent{\bf Proof of Theorem 2.1}.
In fact, in Lemma 4.1 we have shown that under condition $a<a_0\leq a'<b'\leq b_0<b$ and the relationship (2.1), it has
\beqlb\label{sect4-9}\lim\limits_{n\rightarrow+\infty}\frac{-\ln\inf_{x\in[a_0,b_0]} r_{0,n}(a,b,a',b',x,\beta)}{n}=\gamma(b-a,\beta),~~~{\rm a.s.}.\eeqlb
Next, we will divide the proof into four steps.

\

\noindent\emph{Step 1. Showing \beqlb\label{sect4-10}\lim\limits_{n\rightarrow+\infty}\frac{-\ln\sup_{x\in\bfR} r_{0,n}(a,b,a,b,x,\beta)}{n}=\gamma(b-a,\beta),~~~{\rm a.s.}.~~\eeqlb}
In the case of $a<a'\leq a_0<b_0\leq b'<b,$ we have  ~$$\inf_{x\in[a_0,b_0]} r_{0,n}(a,b,a',b',x,\beta)\geq p_{0,n}(a,b,a_0,b_0,x,\beta)$$~ and
$$\inf_{x\in[a_0,b_0]} r_{0,n}(a,b,a',b',x,\beta)\inf_{x\in[a',b']}r_{n,n+1}(a,b,a_0,b_0,x,\beta)\leq p_{0,n+1}(a,b,a_0,b_0,x,\beta).$$
Therefore, we can deduce that if $a<a'\leq a_0<b_0\leq b'<b,$ we also have
\beqlb\label{sect4-11}\lim\limits_{n\rightarrow+\infty}\frac{-\ln\inf_{x\in[a_0,b_0]} r_{0,n}(a,b,a',b',x,\beta)}{n}=\gamma(b-a,\beta).~~~{\rm a.s.}.~~\eeqlb
%$$r_{0,t}(a,b,a,b,x,\beta)\leq \inf\limits_{y\in[x-\epsilon,x+\epsilon]}r_{0,t}(a-\epsilon,b+\epsilon,a-\epsilon,b+\epsilon,y,\beta).$$
Choosing an $\varepsilon>0$ arbitrarily. Let $M:=\lceil\frac{b-a}{\varepsilon}\rceil,$ $y_i=\min\{a+i\varepsilon,b\}, i=0,1,2,\ldots,M.$
\beqlb\label{sect4-12}\sup_{x\in\bfR}r_{0,n}(a,b,a,b,x,\beta)&=&\sup_{x\in[a,b]}r_{0,n}(a,b,a,b,x,\beta)\nonumber
\\&\leq& \max_{0\leq i\leq M-1}\inf_{x\in[y_i,y_{i+1}]}r_{0,n}(a-4\varepsilon,b+4\varepsilon,a-2\varepsilon,b+2\varepsilon,x,\beta)\nonumber
\\&:=&\max_{0\leq i\leq M-1}\overline{r}_{i,n}.\eeqlb
By (\ref{sect4-11}), we know that for each positive integer $i\in[0,M-1]$, it always has
$\liminf_{n\rightarrow+\infty}\frac{-\ln \overline{r}_{i,n}}{n}\geq\gamma(b-a+8\varepsilon, \beta).$ Besides, for fixed $\varepsilon>0,$ $M$ is finite. Thus we have
$$ \liminf\limits_{n\rightarrow+\infty}\frac{-\ln \sup_{x\in\bfR}\bfP^{x}(\forall_{s\in[0,n]}\ \beta W_{s}+a\leq B_{s} \leq \beta W_{s}+b|W)}{n} \geq\gamma(b-a+8\varepsilon,\beta),~~{\rm a.s.}.$$
Moreover, Lemma 4.1 implies that
$$\limsup\limits_{n\rightarrow+\infty}\frac{-\ln \sup_{x\in\bfR}\bfP^{x}(\forall_{s\in[0,n]}\ \beta W_{s}+a\leq B_{s} \leq \beta W_{s}+b|W)}{n} \leq\gamma(b-a,\beta),~~{\rm a.s.}.$$
By Proposition 4.1(1), for each fixed $\beta,$ the function $c\mapsto\gamma(c,\beta)$ is convex hence it is continue. Let $\varepsilon\rightarrow 0,$ we get (\ref{sect4-10}).

\

\noindent\emph{Step 2. Changing time axis from $n\in\bfN$ to $t\in\bfR^{+}.$}

\

Assuming that $t\in(n,n+1).$ Just notice that when $a'<a''<b''<b',$ we have
\beqlb\label{sect4-13}\frac{1}{t}\ln\sup_{x\in\bfR} r_{0,t}(a,b,a,b,x,\beta)\leq \frac{1}{n}\ln\sup_{x\in\bfR} r_{0,n}(a,b,a,b,x,\beta)\eeqlb  and \beqlb\label{sect4-14}\frac{\ln\inf_{x\in[a_0,b_0]}r_{0,t}(a,b,a',b',x,\beta)}{t} &\geq& \frac{\inf_{x\in[a_0,b_0]}r_{0,n}(a,b,a'',b'',x,\beta)}{n+1}\nonumber \\&~~&+\frac{\inf_{x\in[a'',b'']}r_{n,n+1}(a',b',a',b',x,\beta)}{n+1}.\eeqlb

Utilizing (4.3) (4.9) and (4.10), we complete the step {\emph 2}. According to the above discussion, we have shown the almost surely convergence in (2.2). (The only difference is the expression of $\gamma(\beta)$).

\

\noindent\emph{Step 3. Showing the $L^p (p\geq 1)$ convergence in (2.2).}

\

 Because we have proved that the convergence in (2.2) is almost surely,  step {\emph 3} is equivalent to show $\{\frac{-\ln\inf_{x\in[a_0,b_0]}r_{0,n}(a,b,a',b',x,\beta)}{n}\}_{n\in\bfN}$ is $L^p$ uniformly integrable when $a<a_0\leq a'<b'\leq b_0<b$.
Denote $r_{i,i+1}:=-\ln\inf_{x\in[a_0,b_0]}r_{i,i+1}(a,b,a',b',x,\beta)$.
Note that
\beqlb\label{sect4-15}0\leq \frac{-\ln\inf_{x\in[a_0,b_0]}r_{0,n}(a,b,a',b',x,\beta)}{n}\leq \frac{\sum_{i=0}^{n-1}r_{i,i+1}}{n}.\eeqlb Theorem 3.1 shows
$\bfE(r^p_{0,1})<+\infty$ for any $p\geq 1$. Therefore, by the Birkhoff ergodic theorem, we know
$$\frac{\sum_{i=0}^{n-1} r_{i,i+1}}{n}\rightarrow \bfE(r_{0,1}),~~{\rm a.s.~~and ~~in}~~{L^p (p\geq 1)}. $$ Therefore, $\{\frac{\sum_{i=0}^{n-1}r_{i,i+1}}{n}\}_{n\in\bfN}$ is $L^p$ uniformly integrable. By (4.13), we know the sequence $\{\frac{-\ln\inf_{x\in[a_0,b_0]}r_{0,n}(a,b,a',b',x,\beta)}{n}\}_{n\in\bfN}$ is also $L^p$ uniformly integrable. That is to say, we have proved along the discrete time axis $n\in\bfN$ the convergence in (2.2)
is  $L^p ~(p\geq 1).$ Just note the right-hand side of (4.13) and (4.14) are both converge to $\gamma(b-a,\beta)$ in $L^p ~(p\geq 1)$ since (4.3) also holds in the sense of $L^p ~(p\geq 1).$ Hence we can
also change time axis from $n\in\bfN$ to $t\in\bfR^{+}$ by the same way of step {\emph 2}.

\

\noindent\emph{Step 4. Define $\gamma(\beta):=\gamma(1,\beta)$ and show $\gamma(1)\leq 4\pi^2,$ $\gamma(0)=\frac{\pi^2}{2}$ and  $\gamma(\beta)\geq\frac{\pi^{2}(1+\beta^{2})}{2}.$}

\
Firstly, $\gamma(1)\leq 4\pi^2$  can be derived directly from \cite[Corollary 4.4]{DL2005}.

\
According to step \emph{1-2} of this proof, we can see that for any $x\in (a,b),$
$$\lim\limits_{t\rightarrow+\infty}\frac{r_{0,t}(a,b,a,b,x,\beta)}{t}=\gamma(b-a,\beta), ~~~{\rm a.s.}.$$
Moreover, note that for each $t,d>0,$ if $a<0<b,$ we have
$$P^0(\forall_{s\leq t} B_s\in[da+\beta W_s, db+\beta W_s]|W)\stackrel{d}{=}P^0(\forall_{s\leq t} B_{s/d^2}\in[a+\beta W_{s/d^2}, b+\beta W_{s/d^2}|W),$$ where ``$X\stackrel{d}{=}Y$'' means that $X$ and $Y$ have the same distribution. That implies $\gamma(b-a,\beta)=\frac{\gamma(1,\beta)}{(b-a)^2}.$
Therefore, it is reasonable to define $\gamma(\beta):=\gamma(1,\beta).$  So far we have given the whole proof of (2.2).

The only rest thing is to show $\gamma(\beta)\geq\frac{\pi^{2}(1+\beta^{2})}{2}.$ We can use the method which has also been used in the corresponding part in \cite{DL2005}. By the Jensen's inequality we have
$$\bfE(-\ln P^0(\forall_{s\leq t} B_s-\beta W_s\in[-1/2,1/2]|W))> -\ln\bfE(P^0(\forall_{s\leq t} |B_s-\beta W_s|\leq 1/2|W)).$$
%$$\bfE(-\ln P^0(\forall_{s\leq t} B_s-\beta W_s\in[-\frac{1}{2},\frac{1}{2}]|W))> -\ln\bfE(P^0(\forall_{s\leq t} B_s-\beta W_s\in[\frac{1}{2},\frac{1}{2}]|W)),$$
%By step {\emph 3} of this proof, we know
%$$\lim\limits_{t\rightarrow+\infty}\frac{\bfE\big(-\ln P^0(\forall_{s\leq t} B_s\in[a+\beta W_s, b+\beta W_s]|W)\big)}{t}=\frac{\gamma(\beta)}{(b-a)^2}.$$
Let $\tilde{B}$ be a Brownian motion with parameters $\bfE(\tilde{B}_t)=0,\bfE(\tilde{B}^2_t)=(1+\beta^2)t,\forall t\geq 0.$ Then the annealed expectation
$$\bfE(P^0(\forall_{s\leq t} |B_s-\beta W_s|\leq 1/2|W))=\bfP(\forall_{s\leq t} \tilde{B}_s\in[-1/2, 1/2]).$$
It is well known that $$\lim\limits_{t\rightarrow+\infty}\frac{-\ln\bfP(\forall_{s\leq t} |\tilde{B}_s|\leq \frac{1}{2})}{t}=\frac{\pi^2(1+\beta^2)}{2},~~\gamma(0)=\lim\limits_{t\rightarrow+\infty}\frac{-\ln\bfP(\forall_{s\leq t} |B_s|\leq\frac{1}{2})}{t}=\frac{\pi^2}{2}.$$
Hence we have $\gamma(\beta)\geq\frac{\pi^2(1+\beta^2)}{2}.$ Moreover, combining with Proposition 4.1 (2) which shows that $\gamma(\beta)$ is even and convex,
we know $\gamma(\beta)$ is strictly increasing to $+\infty$ on $[0,+\infty)$ and strictly decreasing on $(-\infty,0].$

\qed

\section{Proof of Corollary 2.1 and 2.2}
By scaling property of Brownian motion, it is easy to see that the convergence in (2.4) and (2.5) are in Probability. Thanks to (3.2), we can strengthen it to
almost surely.

\noindent{\bf Proof of Corollary 2.1.} The proof of the upper bound (2.5) is more easier and similar with the lower bound (2.4), so here we only prove (2.4).
We choose an $A>0$ arbitrarily. Denote $M:=\lfloor A^{-1}t^{1-2\alpha}\rfloor, z_i:=iAt^{2\alpha}.$
Without loss of generality, we assume $a_0<a'<b'<b_0$ and choose $a'',b''$ such that $a'<a''<b''<b'.$ It is not hard to see
\beqlb\label{sect5-1} &&\inf_{x\in[a_0t^\alpha,b_0t^\alpha]}\frac{\ln \bfP^x(\forall_{s\leq t} B_s-\beta W_s\in[at^\alpha, bt^\alpha],B_t-\beta W_t\in[a't^\alpha, b't^\alpha]|W)}{t^{1-2\alpha}}\nonumber
\\&\geq&\frac{1}{A}\frac{1}{A^{-1}t^{1-2\alpha}}\Big(\sum_{i=0}^{M-1}V_i(t)+U_{M}(t)\Big),
\eeqlb
where
\beqnn V_i(t)&=&
 \inf_{x\in[a_0t^\alpha,b_0t^\alpha]}\ln\bfP(\forall_{z_i\leq s\leq z_{i+1} } B_s-\beta (W_s-W_{z_i})\in[at^\alpha, bt^\alpha],
 \\  &&B_{z_{i+1}}-\beta (W_{z_{i+1}}-W_{z_i})\in[a''t^\alpha, b''t^\alpha]|W, B_{z_i}=x),\eeqnn
$$U_{M}(t)=\inf_{x\in[a''t^\alpha,b''t^\alpha]}\ln\bfP(\forall_{z_M\leq s\leq z_{M+1}} B_s-\beta (W_s-W_{z_M})\in[a't^\alpha, b't^\alpha]|W,B_{z_M}=x ).$$
Note that for each $t>0$,~$U_{M}(t)\stackrel{d}{=} \inf_{x\in[a'',b'']}\bfP^x(\forall_{s\leq A}B_s-\beta W_s\in[a',b']|W).$  According to Theorem 3.1, we know $\bfE(U^j_{M}(t))<+\infty$  for any $j\in\bfN$, hence $\lim\limits_{t\rightarrow +\infty}\frac{U_{M}(t)}{t^{1-2\alpha}}=0.$
%Define $\lambda:=\bfE(V_{0}(1))=\bfE(\inf_{x\in[a_0,b_0]}\bfP^x(\forall_{s\leq A}B_s-\beta W_s\in[a,b],B_t-\beta W_t\in[a'',b'']|W)).$
We should note that for any fixed $t>0,$ the sequence $V_0(t),V_1(t),...V_{M-1}(t)$ are i.i.d.. Besides, for any fixed $t>0, \forall i\in[0,M-1]\cap\bfN,$ $V_i(t)$ has the same distribution with
$$V_{0}(1)=\inf_{x\in[a_0,b_0]}\bfP^x(\forall_{s\leq A}B_s-\beta W_s\in[a,b],B_A-\beta W_A\in[a'',b'']|W).$$ Now we use Borel-Cantelli 0-1 law to show that $$\lim\limits_{t\rightarrow+\infty}\frac{\sum_{i=0}^{M-1}V_i(t)}{A^{-1}t^{1-2\alpha}}=\bfE(V_{0}(1)).$$

Let $V'_{i}(t):=V_i(t)-\bfE(V_i(t))=V_i(t)-\bfE(V_0(1)).$
Choosing an even positive integer $m$ such that $\frac{(1-2\alpha)m}{2}>1.$ According to (3.2), for any $\varepsilon>0,$ there exists a finite constant $C$ depend on $m$ such that
\beqnn\bfP\Big(\Big|\frac{\sum_{i=0}^{M-1}V'_i(t)}{M}\Big|\geq \varepsilon\Big)&\leq& \bfE\Big(\frac{(\sum_{i=0}^{M-1}V'_i(t))^m}{M^m\varepsilon^{m}}\Big)
\leq \frac{C\mathcal{C}^{m/2}_{M}+o(M^{m/2})}{M^m\varepsilon^{m}},\eeqnn
where $\mathcal{C}$ is the combinatorial number. By Borel-Cantelli 0-1 law, we can obtain $$\lim\limits_{t\rightarrow+\infty}\frac{\sum_{i=0}^{M-1}V_i(t)}{A^{-1}t^{1-2\alpha}}=\bfE(V_0(1)).$$
Combining with (5.1), it implies that for any $A>0,$ we have
\beqnn &&\inf_{x\in[a_0t^\alpha,b_0t^\alpha]}\frac{\ln \bfP^x(\forall_{s\leq t} B_s-\beta W_s \in[at^\alpha, bt^\alpha],B_t-\beta W_t\in[a't^\alpha, b't^\alpha]|W)}{t^{1-2\alpha}}
\\&~~~~~&~~\geq\frac{1}{A}\bfE(\inf_{x\in[a_0,b_0]}\bfP^x(\forall_{s\leq A}B_s-\beta W_s\in[a,b],B_A-\beta W_A\in[a'',b'']|W)),~~~{\rm a.s.}.
\eeqnn
According to the $L^1$ convergence in Theorem 2.1, we get the lower bound (2.4) by taking $A\rightarrow+\infty$.
\qed

%\beqnn
%\leq \frac{\sum_{i=0}^{\lfloor A^{-1}t^{1-2\alpha}\rfloor-1}\sup_{x\in[at^\alpha,bt^\alpha]}\ln\bfP^x(\forall_{iAt^{2\alpha}\leq s\leq (i+1)At^{2\alpha} } B_s-\beta (W_s-W_{it^{2\alpha}})\in[at^\alpha, bt^\alpha]|W, B_{it^{2\alpha}}=x)}{t^{1-2\alpha}/A}\frac{1}{A}
%\\&:=&\frac{\sum_{i=0}^{\lfloor A^{-1}t^{1-2\alpha}\rfloor-1}V_i(t)}{A^{-1}t^{1-2\alpha}}\frac{1}{A}
%\eeqnn

At last, we give the proof of Corollary 2.2.

\noindent{\bf Proof of Corollary 2.2.}
Recalling the notation at the beginning of section 4. The key step of this proof is to observe that for any $u,v>0$,
 it always has
\beqlb\label{sect-2}
 \lim\limits_{t\rightarrow+\infty}\inf\limits_{x\in[a_0,b_0]}\frac{\ln r_{ut,(u+v)t}(a,b,a',b',x,\beta)}{vt}=\frac{-\gamma(\beta)}{(b-a)^2},~~a.s.,~~~~
 \eeqlb
\beqlb\label{sect-3}
 \lim\limits_{t\rightarrow+\infty}\sup_{x\in\bfR}\frac{\ln r_{ut,(u+v)t}(a,b,a,b,x,\beta)}{vt}=\frac{-\gamma(\beta)}{(b-a)^2},~~a.s..~~ \eeqlb
%We will also only show the lower bound (2.6) since the proof of upper bound (2.7) is so analogous and more easier!
Now let us first prove (5.2) and (5.3). For any $m\in\bfN,$ denote $K:=\lfloor\frac{vt}{m}\rfloor,$ $z_k:=ut+km, k\in\bfN.$ Choosing $a'',b''$ such that $a'<a''<b''<b',$ by Markov property
we have
\beqnn
&&\frac{-\ln\inf_{x\in[a_0,b_0]} r_{ut,(u+v)t}(a,b,a',b',x,\beta)}{vt}
\\&\leq& \frac{\sum_{k=0}^{K-1}\overline{G}_k(t)+q_{z_K,z_{K+1}}(a',b',a'',b'',\beta)}{K}\times \frac{K}{vt},
\eeqnn\beqnn
\frac{-\ln\sup_{x\in\bfR} r_{ut,(u+v)t}(a,b,a,b,x,\beta) }{vt}
\geq \frac{\sum_{k=0}^{K-1}\underline{G}_k(t)}{K}\times \frac{K}{vt},
\eeqnn
Where $$\overline{G}_k(t):=-\ln \inf_{x\in[a_0,b_0]} r_{z_k,z_{k+1}}(a,b,a'',b'',x,\beta),$$
$$\underline{G}_k(t):=-\ln \sup_{x\in\bfR} r_{z_k,z_{k+1}}(a,b,a,b,x,\beta).$$
Notice that for any fixed $t>0,$
 $\{\overline{G}_k(t)\}_{k\in \bfN}$ and $\{\underline{G}_k(t)\}_{k\in \bfN}$  are both i.i.d. sequence. And for any $t>0,k\in[0,K-1]\cap\bfN,$
$$\overline{G}_k(t)\stackrel{d}{=}\inf_{x\in[a_0,b_0]} r_{0,m}(a,b,a'',b'',x,\beta),~~\underline{G}_k(t)\stackrel{d}{=}\sup_{x\in\bfR}r_{0,m}(a,b,a,b,x,\beta).$$
 Using Borel-Cantelli 0-1 law, similar with the corresponding part of the proof of Corollary 2.1, we can get (5.2) and (5.3).

Define$$\upsilon:=\min\Big\{a_0-f(0),~g(0)-b_0,~\frac{\inf_{s\in[0,1]}\big(g(s)-f(s)\big)}{3}\Big\}.$$
It is obvious that $\upsilon>0$ since $f(s),g(s)$ are both continue functions in closed interval $[0,1]$ and $f(s)<g(s).$
Moreover,there exists $A_0>0,$  for each $A\geq A_0, A\in \bfN,$ only if $|s_1-s_2|\leq \frac{1}{A},$ we have $$\max\{|f(s_1)-f(s_2)|,|g(s_1)-g(s_2)|\}<\frac{\upsilon}{2}.$$
Denote
$$\overline{Q}_t:=\sup_{x\in\bfR}\bfP^{x}\Big(\forall_{s\leq t}\ \beta W_{s}+f\Big(\frac{s}{t}\Big)\leq B_{s} \leq \beta W_{s}+g\Big(\frac{s}{t}\Big)|W\Big),$$
$$\underline{Q}_t:=\inf\limits_{x\in[a_0,b_0]}\bfP^{x}\Big(\forall_{s\leq t}~f\Big(\frac{s}{t}\Big)\leq B_{s}-\beta W_{s} \leq g\Big(\frac{s}{t}\Big),~a'\leq B_{t}-\beta W_{t}\leq b'|W\Big).$$
%$t_{i,A}:=\frac{it}{A},i=0,1,...A,$\underline{f}_{A,A}=\underline{g}_{A,A}:=a',\overline{f}_{A,A}=\overline{g}_{A,A}:=b'$. And
For $i=0,1,\ldots,A-1,$~define $$\underline{f}_{i,A}:=\inf_{s\in[it/A,(i+1)t/A]}f(s),~~\overline{f}_{i,A}:=\sup_{s\in[it/A,(i+1)t/A]}f(s);$$
$$\underline{g}_{i,A}:=\inf_{s\in[it/A,(i+1)t/A]}g(s),~~\overline{g}_{i,A}:=\sup_{s\in[it/A,(i+1)t/A]}g(s).$$
By Markov property we get
\beqlb\label{sect5-4}\overline{Q}_t\leq \prod_{i=0}^{A-1}\sup_{x\in\bfR}r_{\frac{it}{A},\frac{(i+1)t}{A}}(\underline{f}_{i,A},\overline{g}_{i,A},\underline{f}_{i,A},\overline{g}_{i,A},x,\beta)\eeqlb
and
\beqlb\label{sect5-5}\underline{Q}_t&\geq&\prod_{i=0}^{A-2}\inf_{x\in[f(\frac{i}{A})+\upsilon,g(\frac{i}{A})-\upsilon]}r_{\frac{it}{A},~\frac{(i+1)t}{A}}(\overline{f}_{i,A},\underline{g}_{i,A},
f(\frac{i+1}{A})+\upsilon,g(\frac{i+1}{A})-\upsilon,x,\beta)\nonumber
\\&~~~&\times~\inf_{x\in[f(\frac{A-1}{A})+\upsilon,g(\frac{A-1}{A})-\upsilon]}r_{\frac{(A-1)t}{A},~t}(\overline{f}_{A-1,A},\underline{g}_{A-1,A},
a',b',x,\beta).~~~~~~~~~\eeqlb%:=\prod_{i=1}^{A}\Psi_{i,A}
%-\ln \inf_{x\in[a'',b'']} r_{qt+Km,qt+(K+1)m}(a',b',a',b',x,\beta)
Notice that $$\lim\limits_{A\rightarrow+\infty}\frac{-\sum_{i=0}^{A-1}(\underline{g}_{i,A}-\overline{f}_{i,A})^{-2}}{A}
=\lim\limits_{A\rightarrow+\infty}\frac{-\sum_{i=0}^{A-1}(\overline{g}_{i,A}-\underline{f}_{i,A})^{-2}}{A}
=C_{f,g},$$
Appling (5.2) (5.3) to (5.4) (5.5) we complete the proof of Corollary 2.2.

%$$\geq \prod_{i=1}^{A}\inf_{x\in[f(\frac{i-1}{A})+\upsilon,g(\frac{i-1}{A})-\upsilon]}\bfP^{x}\big(\forall_{\frac{i-1}{A}t\leq s\leq \frac{i}{A}t}B_s-\beta W_s\in[\overline{f}_{i,A},\underline{g}_{i,A}],B_{\frac{i}{A}t}-\beta W_{\frac{i}{A}t}\in[\overline{f}_{i,A}+\upsilon,\underline{g}_{i,A}-\upsilon]\big)$$

%\begin{cor}
%Let $\{a(s)\}$ and $\{b(s)\}$ are two functions from $\bfR^{+}$ to $\bfR^{+}$,   which are both convergence to zero as $s\rightarrow +\infty.$   and
%$a<a_0\leq b_0<b, a\leq a'< b'\leq b,$ Then

 %$$-\lim\limits_{t\rightarrow+\infty}\inf\limits_{x\in[a_0,b_0]}\frac{\ln \bfP(\forall_{s\leq t} B_{s}\in I^{s,\beta}_{a-a(s),b+b(s)},B_{t}\in I^{t,\beta}_{(a',b')}|W)}{t}=\gamma(b-a,\beta)~~~a.s.$$
%\end{cor}
 \ack
I want to thank my supervisor Wenming Hong for his constant concern on my work and giving me a good learning environment. I also want to thank Bastien Mallein
for giving me a lot of valuable advices and useful tips. This research is partly supported by NSFC (NO.11531001, 11626245).


\begin{thebibliography}{99}
\bibitem{CR1979}
M. Cs\"{o}rg\H{o} and P. R\'{e}v\'{e}sz.~How big are the increments of a Wiener process?~
 Acta Mathematica Academiae Scientiarum Hungarica. 33(1-2):37-49, 1979.
\bibitem{D2003}
S. Dereich. Small ball probabilities around random centers of Gaussian measures and applications to quantization. Journal of Theoretical Probability. 16(2):427-449, 2003.
\bibitem{DFMS2003}
S. Dereich, F. Fehringer, A. Matoussi and M. Scheutzow.~ On the link between
small ball probabilities and the quantization problem for Gaussian measures on Banach
spaces. Journal of Theoretical Probability. 16(1):249-265, 2003.
\bibitem{DL2005}
S. Dereich and M. A. Lifshits.~ Probabilities of randomly centered small balls and quantization in Banach spaces. Annals of Probability. 33(4):1397-1421, 2005.
\bibitem{DZ1998}
A. Dembo and O. Zeitouni.~\emph{Large Deviations Techniques and Applications.} Springer-Verlag New York, 1998.
\bibitem{G2002}
R. J. Gardner.~The Brunn-Minkowski inequality. Bulletin of the American Mathematical Society. 39(3):355-405, 2002.
%\bibitem{GHS2011}
%N. Gantert, Y. Hu and Z. Shi. Asymptotics for the survival probability in a killed branching random walk. Ann. Inst. Henri Poincar\'{e} Probab. Stat. 47(1):111-129, 2011.
\bibitem{IM1974}
K. It\^{o} and H. P. McKean Jr.~\emph{Diffusion Processes and Their Sample Paths.} Second printing, corrected, Die Grundlehren der mathematischen Wissenschaften, Band 125. Springer, Berlin, 1974.
\bibitem{K1997}
O. Kallenberg.~\emph{Foundations of Modern Probability.} Probability and its Applications. Springer, New York, 1997.
\bibitem{LS2002}
M. Lifshits and Z. Shi.~The first exit time of Brownian motion from a parabolic domain. Bernoulli. 8(6):745-765, 2002.
\bibitem{L2003}
W. V. Li.~The first exit time of a Brownian motion from an unbounded convex domain. Annals of Probability. 31(2):1078-1096, 2003.
%\bibitem{LSl2011}
%D. Lu, L. Song.~ The first exit time of a Brownian motion
%from the minimum and maximum parabolic domains. Journal
%of Theoretical Probability. 24(4):1028¨C1043, 2011.
\bibitem{LSl2012}
D. Lu and L. Song.~The asymptotic behavior of a Brownian motion with a drift from a random domain. Communications in Statistics. 41(1):62-75, 2012.
\bibitem{Lv2018}
Y. Lv.~Small deviation for random walk with random environment in time. ArXiv e-prints, arXiv:1803.08772, 2018.
\bibitem{MM2015}
B. Mallein and P. Mi{\l}o\'{s}.~Brownian motion and random walks above quenched
random wall. Accepted to Ann. Inst. Henri Poincar\'{e} Probab. Stat. ArXiv e-prints, arXiv:1507.08578, 2015.
\bibitem{N1981}
A. A. Novikov.~On estimates and asymptotic behavior of non-exit probabilities of a Wiener
process to a moving boundary. Mathematics of the USSR-Sbornik. 38(4):539-550, 1981.
%\bibitem{NFK1999}
%Novikov A, Frishling V, Kordzakhia N. Approximations of Boundary Crossing Probabilities for a Brownian Motion[J]. Journal of Applied Probability, 1999, 36(4):1019-1030.
\bibitem{NFK2003}
A. Novikov, V. Frishling and N. Kordzakhia.~Time-dependent barrier options and boundary crossing probabilities. Georgian Mathematical Journal. 10(2):325-334, 2003.
\bibitem{R1977}
L. M. Ricciardi.~\emph{Diffusion Processes and Related Topics in Biology (Lecture Notes in Biomath. 14).}
 Springer, Berlin, 1977.
 \bibitem{S1985}
D. Siegmund.~\emph{Sequential Analysis:~Tests and Confidence Intervals.} Springer, New York, 1985.
\bibitem{XZ2015}
L. Xu and D. Zhu.~On the distribution of first exit time for Brownian motion with double linear time-dependent barriers. Isrn Applied Mathematics. 2013(1):64-68, 2013.
\end{thebibliography}
\end{document}